\documentclass{article}

\usepackage{microtype}
\usepackage{graphicx}
\usepackage{subfigure}
\usepackage{booktabs} 
\usepackage[usestackEOL]{stackengine}

\usepackage{hyperref}

\newcommand{\R}{\mathbb{R}}

\def\<#1,#2>{\langle #1,#2\rangle}

\newcommand{\norm}[1]{\|#1\|}
\newcommand{\sqn}[1]{\norm{#1}^2}

\newcommand{\cR}{\mathcal{R}}

\newcommand{\cO}{\mathcal{O}}

\usepackage[accepted]{icml2023}

\usepackage{amsmath}
\usepackage{amssymb}
\usepackage{mathtools}
\usepackage{amsthm}

\usepackage[capitalize,noabbrev]{cleveref}

\theoremstyle{plain}
\newtheorem{theorem}{Theorem}[section]

\newtheorem{lemma}[theorem]{Lemma}
\newtheorem{corollary}[theorem]{Corollary}
\theoremstyle{definition}

\newtheorem{assumption}[theorem]{Assumption}
\theoremstyle{remark}

\usepackage[textsize=tiny]{todonotes}

\newcommand{\sX}{\R^{d_x}}
\newcommand{\sY}{\R^{d_y}}
\newcommand{\ol}[1]{\overline{#1}}
\newcommand{\ul}[1]{\underline{#1}}
\DeclareMathOperator*{\argmin}{arg\,min}

\icmltitlerunning{Submission and Formatting Instructions for ICML 2023}

\begin{document}

\twocolumn[
\icmltitle{Optimal algorithm for min-min strongly convex optimization}



\icmlsetsymbol{equal}{*}

\begin{icmlauthorlist}
\icmlauthor{Dmitry Kovalev}{yyy}
\icmlauthor{Alexander Gasnikov}{comp,comp1,comp2}
\icmlauthor{Grigory Malinovsky}{sch}
\end{icmlauthorlist}

\icmlaffiliation{yyy}{Université catholique de Louvain (UCL), CORE(IMMAQ),}
\icmlaffiliation{comp}{Moscow Institute of Physics and Technology,}
\icmlaffiliation{comp1}{ISP RAS Research Center for Trusted Artificial Intelligence,}
\icmlaffiliation{comp2}{HSE University,}
\icmlaffiliation{sch}{King Abdullah University of Science and Technology}

\icmlcorrespondingauthor{Dmitry Kovalev}{dakovalev1@gmail.com}
\icmlcorrespondingauthor{Alexander Gasnikov}{gasnikov@yandex.ru}
\icmlcorrespondingauthor{Grigory Malinovsky}{grigorii.malinovskii@kaust.edu.sa}
\icmlkeywords{Machine Learning}

\vskip 0.3in
]



\printAffiliationsAndNotice{\icmlEqualContribution} 

\begin{abstract}
Assume that we have to minimize $f(x,y)$, where $f$ is smooth and $\mu_x$-strongly convex in $x$ and $\mu_y$-strongly convex in $y$. 
Optimal fast gradient method of Yurii Nesterov requires $\sim \left(\min\left(\mu_x,\mu_y\right)\right)^{-1/2}$ calculations of $\nabla_x f$ and $\nabla_y f$. 
In this paper we propose an algorithm that requires $\sim \mu_x^{-1/2}$ calculations of $\nabla_x f$  and $\sim \mu_y^{-1/2}$ calculations of $\nabla_y f$.  
In some applications $\mu_x \gg \mu_y$ and calculation of $\nabla_y f$ is much cheaper than $\nabla_x f$. In this case, we have significant wall-clock time acceleration. 
\end{abstract}

\section{Introduction}
Optimal ("black-box") algorithms for basic classes of convex optimization problems were developed many years ago \cite{nemirovski1983problem}. Modern results assume the additional structure of the problem "to look inside the black box" \cite{nesterov2018lectures}. Many results of this type relate to the composite structure of the problem\footnote[2]{Function $F$ is $\mu$-strongly convex, $f$ and $g$ have Lipschitz continuous gradients with constants correspondingly  $L_f$ and $L_g$.}
$$\min_x F(x):= f(x) + g(x),$$
where it is possible to "split" the complexity as $\sim \sqrt{L_f/\mu}$ calculations of $\nabla f$ and  $\sim \sqrt{L_g/\mu}$ calculations of $\nabla g$ \cite{lan2016gradient,ivanova2022oracle,kovalev2022optimal}. But there are lack of "optimal results" for \textit{min-min} problem   
\begin{equation}\label{eq:based_problem}
\min_{x,y} f(x,y),
\end{equation}
where smoothness constants, strong convexity constants, complexities of $\nabla_x f$ and $\nabla_y f$ can be quite different for different blocks $x$ and $y$ as well as the dimensions of these blocks.
For example, such problems play an important role in modeling combined trip distribution, and assignment \cite{de2005solving,gasnikov2014three} in soft clustering \cite{nesterov2020soft}. The natural application in Machine Learning we demonstrate following by Yahoo! Click-prediction model from \cite{dvurechensky2022hyperfast}
\begin{align*}
    \min_{x\in \R^d} f(x)&: = \frac{1}{m}\sum_{k=1}^m \log\left(1 + \exp\left(-\eta^k\langle\xi^k,x\rangle\right)\right)\\
    &+ \lambda_S \sum_{i \in I_S} x_i^2 + \lambda_D \sum_{i \in I_D} x_i^2,
\end{align*}

where $I_S\cup I_D = \left\{1,...,d \right\}$, $I_S\cap I_D = \varnothing$, $|I_S|\gg |I_S|$ and $\lambda_S \gg \lambda_D$.
Here it is natural to put $x:=\left\{x_i\right\}_{i\in I_S}$ and $y:=\left\{x_i\right\}_{i\in I_D}$ in representation \eqref{eq:based_problem}. 

\subsection{Problem formulation and informal formulation of the main result}
Formally speaking, we will consider in this paper the following class of problems
\begin{equation}\label{eq:main}
	\min_{x \in \R^{d_x},y \in \R^{d_y}} f(x,y),
\end{equation}
where $f(x,y)\colon \sX\times \sY \rightarrow \R$ is a convex function that satisfies the Assumptions below.

\begin{assumption}\label{ass:L}
	Function $f(x,y)$ is $(L_x, L_y)$-smooth with $L_x,L_y > 0$. That is, for all $x_1,x_2 \in \sX$, $y_1,y_2 \in \sY$ the following inequality holds:
	\begin{equation}
		\begin{split}
				f(x_2,y_2) &\leq f(x_1,y_1) + \<\nabla_x f(x_1,y_1),x_2-x_1>\\&  \quad+\<\nabla_y f(x_1,y_1),y_2-y_1> \\&\quad+ \frac{L_x}{2}\sqn{x_2-x_1} + \frac{L_y}{2}\sqn{y_2-y_1}.
		\end{split}
	\end{equation}
	
\end{assumption}

Note, that up to a multiplayer $2$ we can consider $L_x$ to be a Lipschitz gradient constant in block $x$ and $L_y$ to be a Lipschitz gradient constant in block $y$.

\begin{assumption}\label{ass:mu}
	Function $f(x,y)$ is $(\mu_x, \mu_y)$-strongly convex. That is, for all $x_1,x_2 \in \sX$, $y_1,y_2 \in \sY$ the following inequality holds:
	\begin{equation}
		\begin{split}
			f(x_2,y_2) &\geq f(x_1,y_1) + \<\nabla_x f(x_1,y_1),x_2-x_1>\\
   &\quad+\<\nabla_y f(x_1,y_1),y_2-y_1> \\&\quad+ \frac{\mu_x}{2}\sqn{x_2-x_1} + \frac{\mu_y}{2}\sqn{y_2-y_1}.
		\end{split}
	\end{equation}
	
\end{assumption}

Informally the main result of the paper is a description of the Block Accelerated Method (BAM, see Section~\ref{sec:Algorithm}) that finds a solution of \eqref{eq:main} with relative precision $\epsilon$ with:
\begin{center}
\framebox{\Longstack[c]{$\cO\left(\sqrt{\frac{L_x}{\mu_x}}\log \frac{1}{\epsilon}\right)$ calculations of $\nabla_x f$\\\\ and $\cO\left(\sqrt{\frac{L_y}{\mu_y}}\log \frac{1}{\epsilon}\right)$ calculations of $\nabla_y f$.}}
\end{center}
According to \cite{nemirovski1983problem,nesterov2018lectures}, it is pretty evident that these bounds correspond to the lower bound separately and even more so together. 

By using the regularization trick (e.g. see \cite{gasnikov2016efficient}), we can reduce the convex problem, but not strongly convex in some block(s) case(s) (denoted by $\circ$), to strongly convex one(s) with $\mu_{\circ} \sim \epsilon/R^2$, where $R$ is a distance in $2$-norm between starting point in the closest solution.

\subsection{Related works}
The close problem formulation considered when studied accelerated coordinate-descent methods \cite{nesterov2012efficiency,richtarik2014iteration,nesterov2017efficiency,ivanova2021adaptive}.\footnote{Note that the the first accelerated gradient schemes \cite{nesterov2012efficiency,richtarik2014iteration} do not allow to obtain this result. The first time this results was announced in \cite{nesterov2015polyak80} by using special coordinate-wise randomization: $p_x \sim \sqrt{L_x}$ and $p_y \sim \sqrt{L_y}$. Next it was developed at different works with different variations \cite{gasnikov2015accelerated,allen2016even,nesterov2017efficiency}.}

However, the result is quite different: with probability $\ge 1 - \delta$
\begin{center}
$\cO\left(\sqrt{\frac{{L}_x}{\min\left\{\mu_x,\mu_y\right\} } }\log \frac{1}{\epsilon}\log \frac{1}{\delta}\right)$ calculations of $\nabla_x f$  and $\cO\left(\sqrt{\frac{{L}_y}{\min\left\{\mu_x,\mu_y\right\}}}\log \frac{1}{\epsilon}\log \frac{1}{\delta}\right)$ calculations of $\nabla_y f$,
\end{center}
where ${L}_x$ is Lipschitz constant of $\nabla_x f(x,y)$ as a function of $x$ and ${L}_y$ is Lipschitz constant of $\nabla_y f(x,y)$ as a function of $y$. 
The same results, but with worse smoothness constants, hold for accelerated alternating methods \cite{beck2017first,diakonikolas2018alternating,guminov2021combination,tupitsa2021alternating}.

Note that by using re-scaling of $y$ variables $y':=\sqrt{\mu_y/\mu_x} y$ it is possible to equalise the constants of strong convexity $\mu_x = \mu_{y'}$. Applying accelerated coordinate-descent from \cite{nesterov2017efficiency} to the re-scaled problem we can obtain in standard variables: 
\begin{center}
$\cO\left(\sqrt{\frac{{L}_x}{\mu_x} }\log \frac{1}{\epsilon}\log \frac{1}{\delta}\right)$ calculations of $\nabla_x f$  and $\cO\left(\sqrt{\frac{{L}_y}{\mu_y}}\log \frac{1}{\epsilon}\log \frac{1}{\delta}\right)$ calculations of $\nabla_y f$.
\end{center}
This result is close to our result, but our approach is deterministic that is better in $\log \frac{1}{\delta}$ times. Moreover, our results are based on very different ideas.  

In the cycle of works  
\cite{bolte2020ah,gladin2021solvingMinMin,gladin2021solving,ostroukhov2022tensor} it was proposed to solve \eqref{eq:main} an outer optimization problem in $x$ block with inexact gradient oracle determined by the solution of the inner problem in $y$ block:
\begin{align*}
    \min_x F(x):&= \min_y f(x,y),\\
    \nabla F(x) &= \nabla_x f(x,y(x))\\
    &= \frac{\partial f}{\partial x} \left(x,y\right)\Big|_{y=y(x)},
\end{align*}
    where $y(x)$ is determined as the solution of $\min_y f(x,y)$.

The most practical results were obtained in the case when $x \in Q \subset \R^{d_x}$, where $d_x$ is small: 
\begin{center}
$\tilde{\cO}\left(d_x\log \frac{1}{\epsilon}\right)$ calculations of $\nabla_x f$  and $\tilde{\cO}\left(d_x\sqrt{\frac{L_y}{\min\left\{\mu_x,\mu_y\right\}}}\log^2 \frac{1}{\epsilon}\right)$ calculations of $\nabla_y f$.
\end{center}
Note that the known lower bound assumes: 
\begin{center}
$\cO\left(d_x\log \frac{1}{\epsilon}\right)$ calculations of $\nabla_x f$  and $\cO\left(\sqrt{\frac{L_y}{\min\left\{\mu_x,\mu_y\right\}}}\log \frac{1}{\epsilon}\right)$ calculations of $\nabla_y f$.
\end{center}
This lower bound may not be tight. 

These results were further generalized to different types of oracles for the inner problem (gradient-free \cite{gladin2021solving}, randomized variance-reduced \cite{gladin2021solvingMinMin}, tensor \cite{ostroukhov2022tensor}). Nevertheless, if $d_x$ is not small, the outer method should be an accelerated one, and such an approach gives:  
\begin{center}
  $\cO\left(\sqrt{\frac{L_x}{\mu_x}}\log \frac{1}{\epsilon}\right)$ calculations of $\nabla_x f$  and $\cO\left(\sqrt{\frac{L_x L_y}{\mu_x\mu_y}}\log^2 \frac{1}{\epsilon}\right)$ calculations of $\nabla_y f$.
\end{center}
That is much worse for the required number of $\nabla_y f$ calculations. So before our work, we knew nothing about the possibility of full-fledged deterministic splitting of complexities into blocks.

\begin{algorithm*}
	\caption{\textbf{B}lock \textbf{A}ccelerated \textbf{M}ethod (BAM)}
	\begin{algorithmic}\label{alg:sliding}
		\STATE {\bf Parameters:} $\eta_x,\eta_y > 0$, $\theta_x, \theta_y > 0$, $\alpha \in (0,1)$ 
		\STATE {\bf Input:} $x^0 = \ol{x}^0 \in \sX$, $y^0 = \ol{y}^0$
		\FOR{$k=0,1,\ldots, K-1$}
		\STATE $\ul{x}^k = \alpha x^k + (1-\alpha)\ol{x}^k $
		\STATE $\ul{y}^k = \alpha y^k + (1-\alpha)\ol{y}^k$
		\STATE find $\ol{y}^{k+1}$ such that
		\begin{equation}\label{eq:ms}
			\norm{\nabla_y f(\ul x^k, \ol y^{k+1}) + (\eta_y\alpha)^{-1}(\ol y^{k+1}-\ul y^k)} \leq (\eta_y\alpha)^{-1}\norm{\ol y^{k+1} - \ul y^k}.
		\end{equation}
		\STATE $\ol{x}^{k+1} = \ul x^k - \eta_x\alpha \nabla_x f(\ul x^k, \ol y^{k+1})$
		\STATE $x^{k+1} = x^k + \alpha(\ul{x}^k - x^{k+1}) - \eta_x \nabla_x f(\ul x^k, \ol y^{k+1})$
		\STATE $y^{k+1} = y^k + \alpha(\ol{y}^{k+1} - y^{k+1}) - \eta_y \nabla_y f(\ul x^k, \ol y^{k+1})$
		\ENDFOR
	\end{algorithmic}
\end{algorithm*}
\section{Main Algorithm}\label{sec:Algorithm}
The BAM algorithm was inspired by a series of recent NeurIPS 2022 papers \cite{kovalev2022optimal,kovalev2022first,kovalev2022first_high} (see also \cite{ivanova2021adaptive,gasnikov2021accelerated,carmon2022recapp}), where the authors use inner-loop (catalyst-type) to obtain optimal accelerated methods for saddle-point problems and high-order methods. We emphasize that BAM is a significantly different algorithm since we need to split the complexities in block arguments.

\begin{lemma}\label{lem:descent}
	Let $\eta_x$ satisfy
	$
		\eta_x \leq (\alpha L_x)^{-1}.
	$
	Then, the following inequality holds:
	\begin{align}
	    		-f(\ul x^{k}, \ol y^{k+1}) &\leq - f(\ol x^{k+1}, \ol y^{k+1})\\\notag
       &\quad-\frac{\eta_x\alpha}{2} \sqn{\nabla_x f(\ul x^{k}, \ol y^{k+1})}.
	\end{align}

\end{lemma}

\begin{theorem}\label{thm:sliding}
	Let $\cR_x^k = \sqn{x^k - x^*}$, $\cR_y^k = \sqn{y^k - y^*}$.
	Let $\Psi^k$ be the following Lyapunov function:
	\begin{align}
		\Psi^k &= (1+\alpha)\left(\eta_x^{-1}\cR_x^{k} + \eta_y^{-1}\cR_y^{k}\right)\\\notag
		&\quad+
		\frac{2}{\alpha}\left( f(\ol x^{k}, \ol y^{k}) - f(x^*,y^*)\right).
	\end{align}
	Let parameters $\eta_x,\eta_y,\alpha$ be defined as follows:
	\begin{equation}
		\alpha = \sqrt{\frac{\mu_x}{L_x}}, \qquad \eta_x = \frac{1}{\sqrt{\mu_x  L_x}}, \qquad \eta_y =\frac{1}{\mu_y}\sqrt{ \frac{\mu_x}{L_x}}.
	\end{equation} 
	Then, iterations of \Cref{alg:sliding} satisfy the following inequality:
	\begin{equation}
		\Psi^{k+1} \leq (1+\alpha)^{-1}\Psi^k.
	\end{equation}
\end{theorem}

\section{Inner Algorithm}

Let us define the following auxiliary function $A^k(y)\colon \sY \rightarrow \R$:
\begin{equation}
	A^k(y) = f(\ul x^k, y) + \frac{1}{2\eta_y\alpha}\sqn{y - \ul y^k}.
\end{equation}
Then, condition \eqref{eq:ms} in \Cref{alg:sliding} can be equivalently written as follows
\begin{equation}\label{eq:ms2}
	\norm{\nabla A^k (\ol y^{k+1})} \leq (\eta_y\alpha)^{-1}\norm{\ol y^{k+1} - \ul y^k}.
\end{equation}
To find $\ol{y}^{k+1}$ that satisfies this condition, we will apply an optimal algorithm for gradient norm reduction \citep{diakonikolas2022potential,kim2021optimizing} to the following minimization problem:
\begin{equation}\label{eq:aux}
	\min_{y \in \sY} A^k(y).
\end{equation}

The following theorem is taken from Remark 1 of \cite{nesterov2021primal}.
\begin{theorem}
	There exists a certain algorithm that, when applied to problem~\eqref{eq:aux} with the starting point $\ul y^k$, returns $\ol y^{k+1}$ satisfying 
	\begin{equation}\label{eq:inner}
		\norm{\nabla A^k(\ol y^{k+1})} \leq \frac{C\max\{L_y,( \eta_y\alpha)^{-1}\}\norm{\ul y^k- \ul y^\star}}{T^2},
	\end{equation}
	where  $T$ is the number of $\nabla A^k(y)$ calls by the algorithm, $\ul{ y}^\star = \argmin_{y \in \sY} A^k(y)$ and $C>0$ is a universal constant.
\end{theorem}

\begin{corollary}\label{cor:inner}
	To output $\ol y^{k+1}$ satisfying condition~\eqref{eq:ms2}, inner algorithm requires the following number of iterations:
	\begin{equation}\label{eq:T}
		T= \sqrt{2C}\max\{1, \sqrt{\eta_y\alpha L_y}\}
	\end{equation}
\end{corollary}
There is a simple approach to achieve the optimal rate $\left(\mathcal{O}\left(\frac{1}{T^2}\right)\right)$ for gradient norm reduction under the initial distance condition. The algorithm runs Nesterov Accelerated Gradient for the first $N/2$ iterations and then it runs OGM-G (Algorithm~\ref{alg:OGM-G}) for the second $N/2$ iterations.  

OGM-G algorithm uses a triangular matrix $\Tilde{\theta}_i$, which determines coefficients for iterations. The first line of the algorithm is the gradient step and the second line is acceleration step using previous points and coefficients $\Tilde{\theta}_i$. 

\begin{algorithm*}
	\caption{\textbf{O}ptimized \textbf{G}radient \textbf{M}ethod (OGM-G)}
	\begin{algorithmic}\label{alg:OGM-G}
		\STATE {\bf Parameters:} stepsize $\gamma$, matrix $\Tilde{\theta}_i$:
 \STATE $\tilde{\theta}_i= \begin{cases}\frac{1+\sqrt{1+8 \tilde{\theta}_{i+1}^2}}{2}, & i=0, \\ \frac{1+\sqrt{1+4 \tilde{\theta}_{i+1}^2}}{2}, & i=1, \ldots, N-1, \\ 1, & i=N,\end{cases}$
		\STATE {\bf Input:} $x^0 = y^0 \in \mathbb{R}^d$
		\FOR{$i=0,1,\ldots, N-1$}
			\STATE $y_{i+1}=x_i-\gamma\nabla f\left(x_i\right) $
			\STATE $x_{i+1}=y_{i+1}+\frac{\left(\tilde{\theta}_i-1\right)\left(2 \tilde{\theta}_{i+1}-1\right)}{\tilde{\theta}_i\left(2 \tilde{\theta}_i-1\right)}\left(y_{i+1}-y_i\right)+\frac{2 \tilde{\theta}_{i+1}-1}{2 \tilde{\theta}_i-1}\left(y_{i+1}-x_i\right)$
		\ENDFOR
	\end{algorithmic}
\end{algorithm*}

\section{Total Complexity}

\Cref{thm:sliding} implies that to find an $\epsilon$-accurate solution of problem~\eqref{eq:main}, \Cref{alg:sliding} requires the following number of calls of $\nabla_x f(x,y)$:
\begin{equation}
	K = \cO\left(\sqrt{\frac{L_x}{\mu_x}}\log \frac{1}{\epsilon}\right).
\end{equation}
\Cref{cor:inner} with the choice of parameters of \Cref{alg:sliding} from \Cref{thm:sliding} implies that the number of inner iterations is the following
\begin{equation}
    \begin{split}
        	T &= \cO\left(\max\{1, \sqrt{\eta_y\alpha L_y}\}\right)\\
 &=\cO\left(\max\left\{1, \sqrt{\frac{L_y\mu_x}{L_x\mu_y}}\right\}\right).
    \end{split}
\end{equation}

Hence, the total number of calls of $\nabla_y f(x,y)$:
\begin{align*}
	K\times T
	&=
	\cO\left(\sqrt{\frac{L_x}{\mu_x}}\log \frac{1}{\epsilon}\right) \times  \cO\left(\max\left\{1, \sqrt{\frac{L_y\mu_x}{L_x\mu_y}}\right\}\right)
	\\&=\cO\left(\max\left\{\sqrt{\frac{L_x}{\mu_x}},\sqrt{\frac{L_y}{\mu_y}}\right\}\log \frac{1}{\epsilon}\right).
\end{align*}
\section{Federated Learning application}
\subsection{Collaborative learning}
Federated learning is a general machine learning framework in which several clients (workers) train the model in a distributed setting without sharing clients' data to maintain privacy~\cite{mcmahan2017communication}. Typically, the training data are distributed across many clients, and these clients can communicate only to the central server (centralized regime)~\cite{konevcny2016federated}. In contrast, in a decentralized regime~\cite{koloskova2020unified}, clients communicate according to a communication graph without any central node. For example, clients can jointly train a text prediction model for a keyboard application without sharing the local sensitive data with other clients or servers. Federated Learning is deployed and implemented in production in cross-device settings~\cite{hard2018federated} and cross-silo settings~\cite{rieke2020future}.

The standard approach is to train a single global model using local updates from clients. The most popular baseline is FedAvg~\cite{khaled2020tighter,woodworth2020local}. The idea of this method is to utilize several local gradient steps before aggregation to reduce communication cost, which is known to be the bottleneck in this scheme. However, due to data heterogeneity, FedAvg has poor convergence guarantees if additional assumptions about similarity are not supposed. In order to correct this issue, several methods were introduced~\cite{karimireddy2020scaffold,mitra2021linear,gorbunov2021local}, and this family of methods has a linear convergence rate in a deterministic regime. However, the communication complexities of such methods are not better than the complexity of the vanilla GD method due to the small stepsizes appeared in the analysis. In recent work~\cite{mishchenko2022proxskip}, it is shown that local steps lead to communication acceleration and subsequent works apply such mechanism to different settings~\cite{malinovsky2022variance,grudzien2022can,condat2022provably}.
\begin{figure*}[t!]
	\centering
	\begin{tabular}{ccc}
		\includegraphics[width=0.31\linewidth]{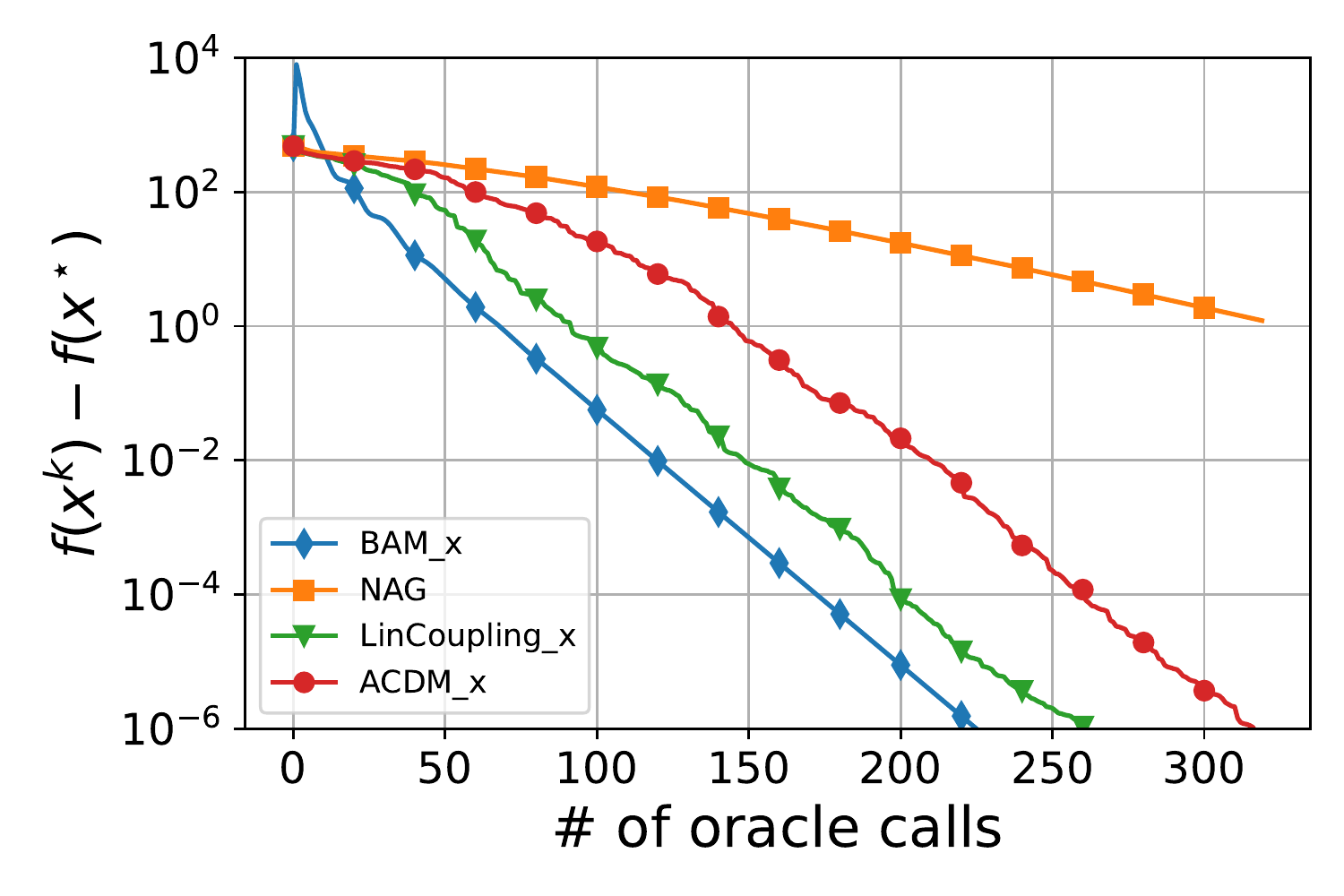}&	\includegraphics[width=0.31\linewidth]{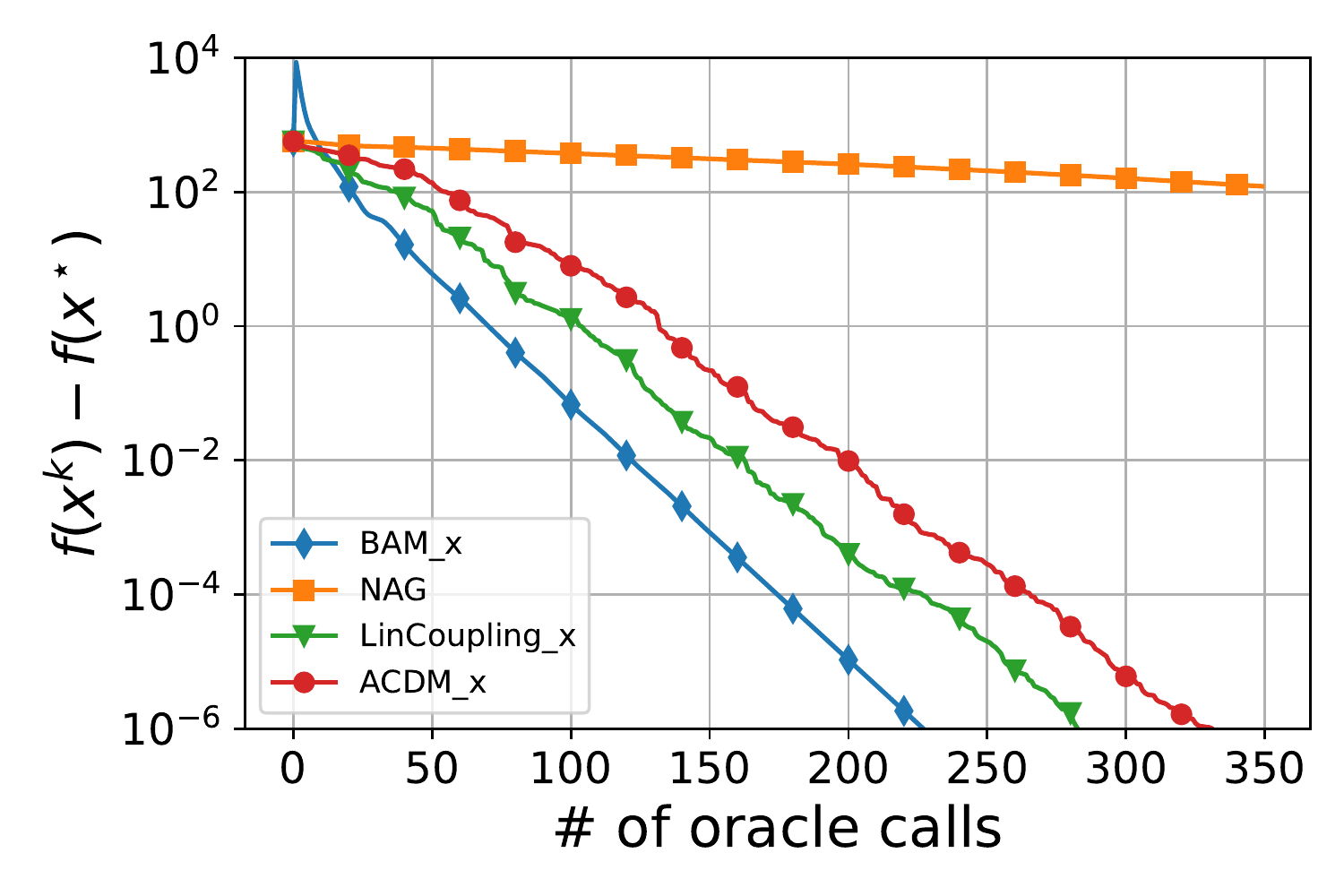}&
		\includegraphics[width=0.31\linewidth]{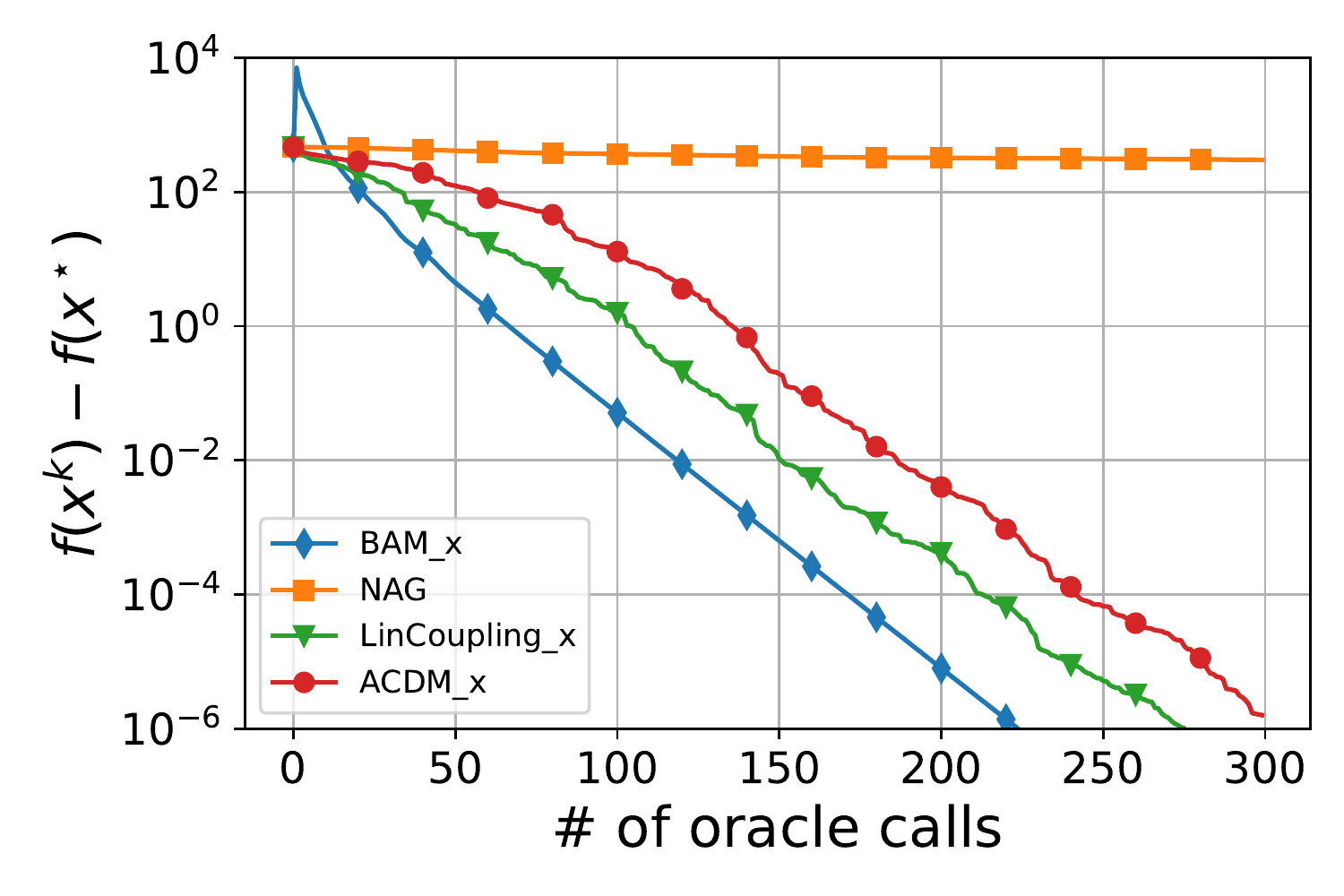}\\
  		\includegraphics[width=0.31\linewidth]{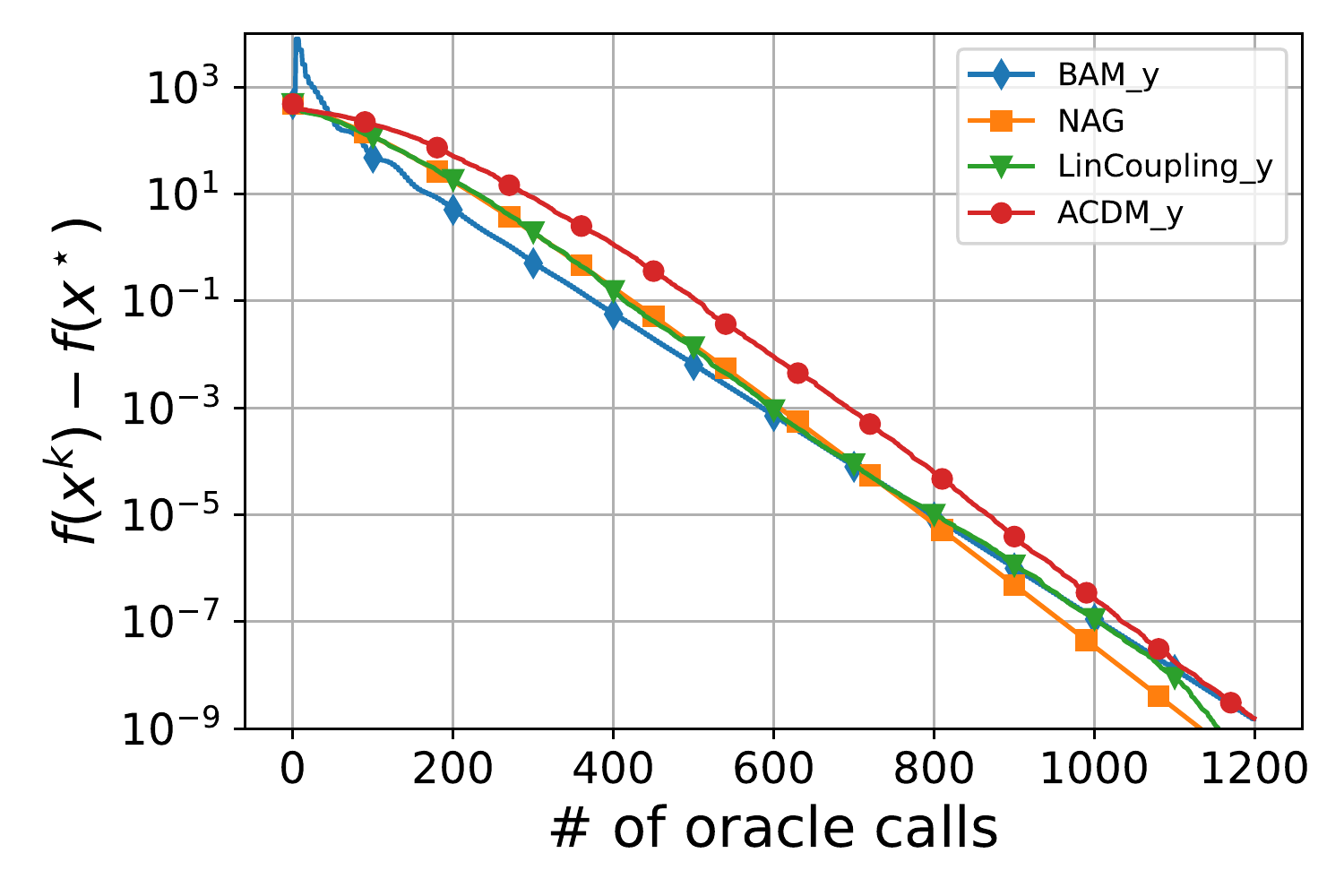}&	\includegraphics[width=0.31\linewidth]{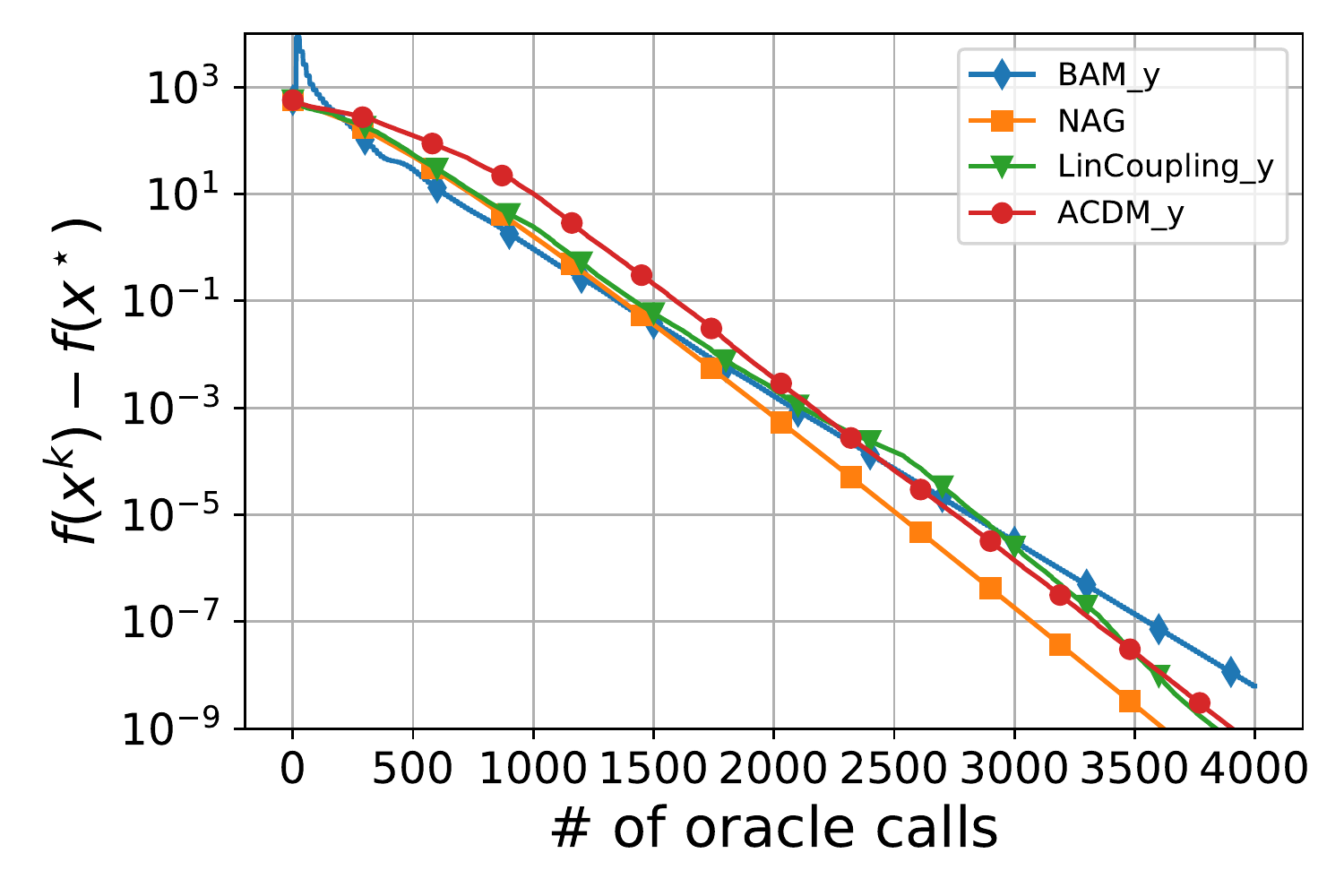}&
		\includegraphics[width=0.31\linewidth]{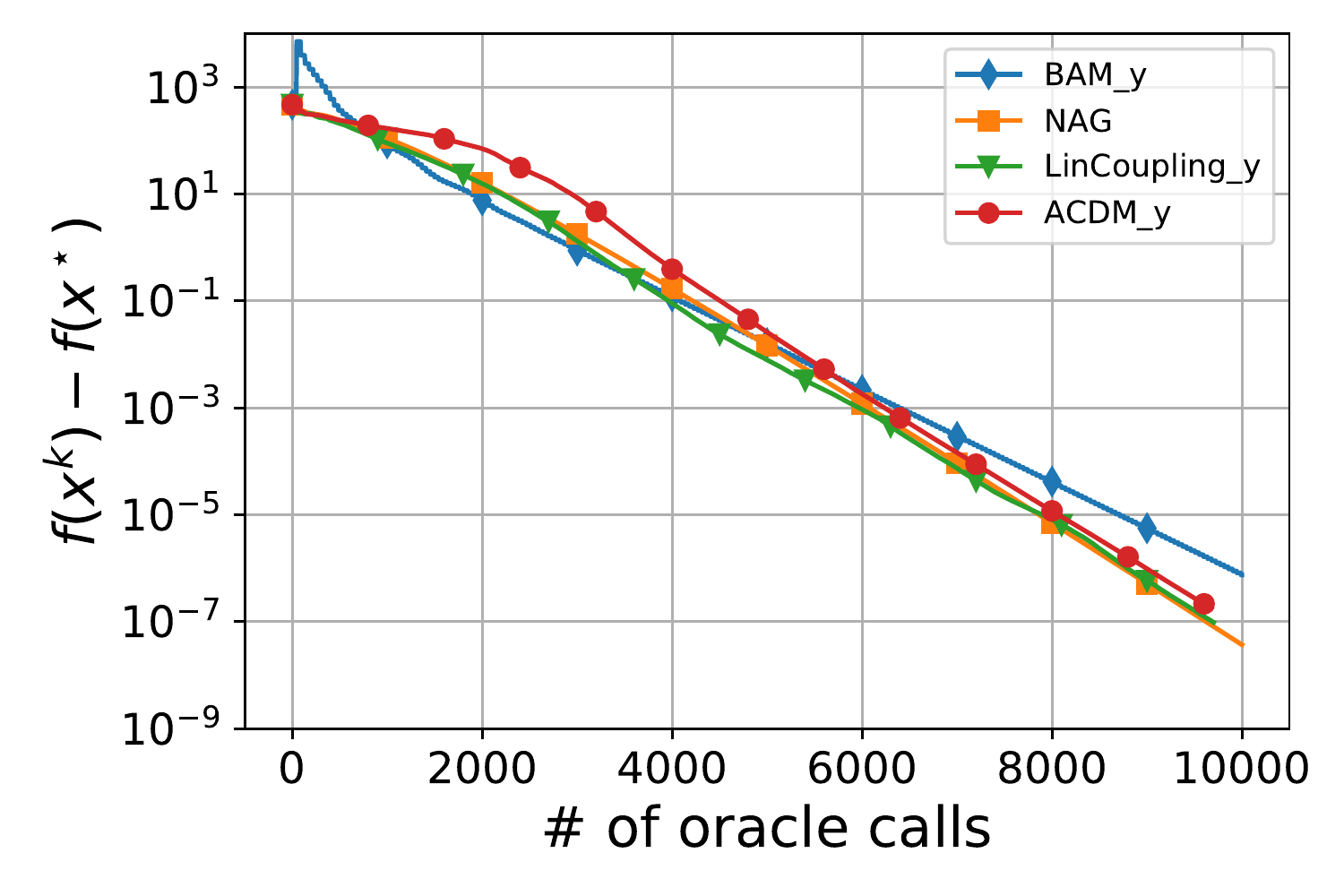}
	\end{tabular}
	\centering
	\caption{Comparison of Block Accelerated Method (BAM), Nesterov Accelerated Method (NAG), Accelerated Coordinate Descent Method (ACDM) and Linear Coupling method (LinCoupling) on quadratic functions. First line represents rate in terms of the $\nabla_x f(x,y)$ oracle calls and the second one represents rate in terms of the $\nabla_y f(x,y)$ oracle calls. We set $L_y = 500$ (left column), $L_y = 5000$ (middle column) and $L_y = 50000$ (right column).  }
	\label{ris:image1}
\end{figure*}

However, global model training can be prohibited in some settings even without sharing data due to privacy constraints. For example, using client-specific embeddings can reveal user identity, which is not allowed by a privacy policy. In order to fix this issue, a concept of partial federated learning was introduced~\cite{singhal2021federated}. In this approach, models have two blocks of parameters: global block $x$ and local blocks $y_i$, which never leave the clients. This technique enables to have interpolation between distributed and non-distributed learning. Partial federated learning is closely connected to personalizing and meta-learning algorithms. The most popular meta-learning algorithm is MAML~\cite{finn2017model}, and connection to federated learning was established in several works~\cite{nichol2018first, chen2018federated, fallah2020personalized}. 

\subsection{Federated Reconstruction}
Let us describe the baseline of partial federated learning called Federated Reconstruction. We have two blocks of coordinates in this framework: user-specific parameters $y_i$ and non-user-specific $x$ parameters. For every communication round, $t$ server sends the global part of parameters $x_t$ to all clients, and then each client reconstructs local parameters $y^i_t$ using the current global model $x_t$. The reconstruction process usually requires several steps. Once the local model is restored, each client updates its copy of global parameters and then sends only updated copies of the global model to the server. The server aggregates these updates and forms the next iterate $x_{t+1}$. 

The new BAM algorithm can be generalized to minimize $f(x,y_1,\ldots,y_M)$ in a distributed setting, and this method can be applied to Federated Reconstruction. Since the communication complexity is proportional to the number of calls of $\nabla_x f(x,y_1,\ldots,y_M)$, the communication complexity is $\mathcal{O}\left( \sqrt{\frac{L_x}{\mu_x}} \log \frac{1}{\varepsilon}\right)$. The communication bottleneck can be overcome in case of a small condition number for local parameters. Moreover, this communication complexity will be optimal. 

\begin{figure*}[t!]
	\centering
	\begin{tabular}{ccc}
        \includegraphics[width=0.31\linewidth]{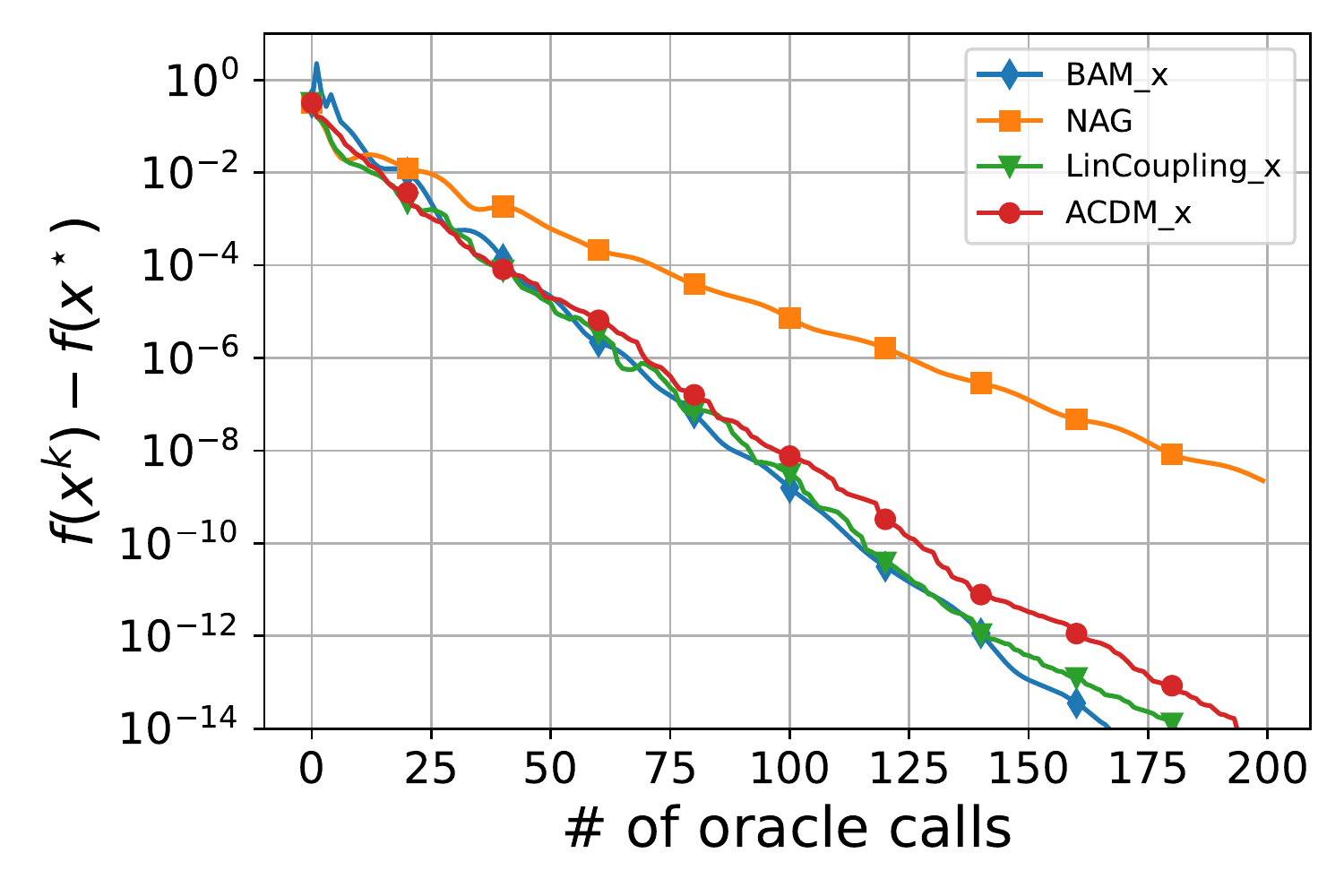}&
		\includegraphics[width=0.31\linewidth]{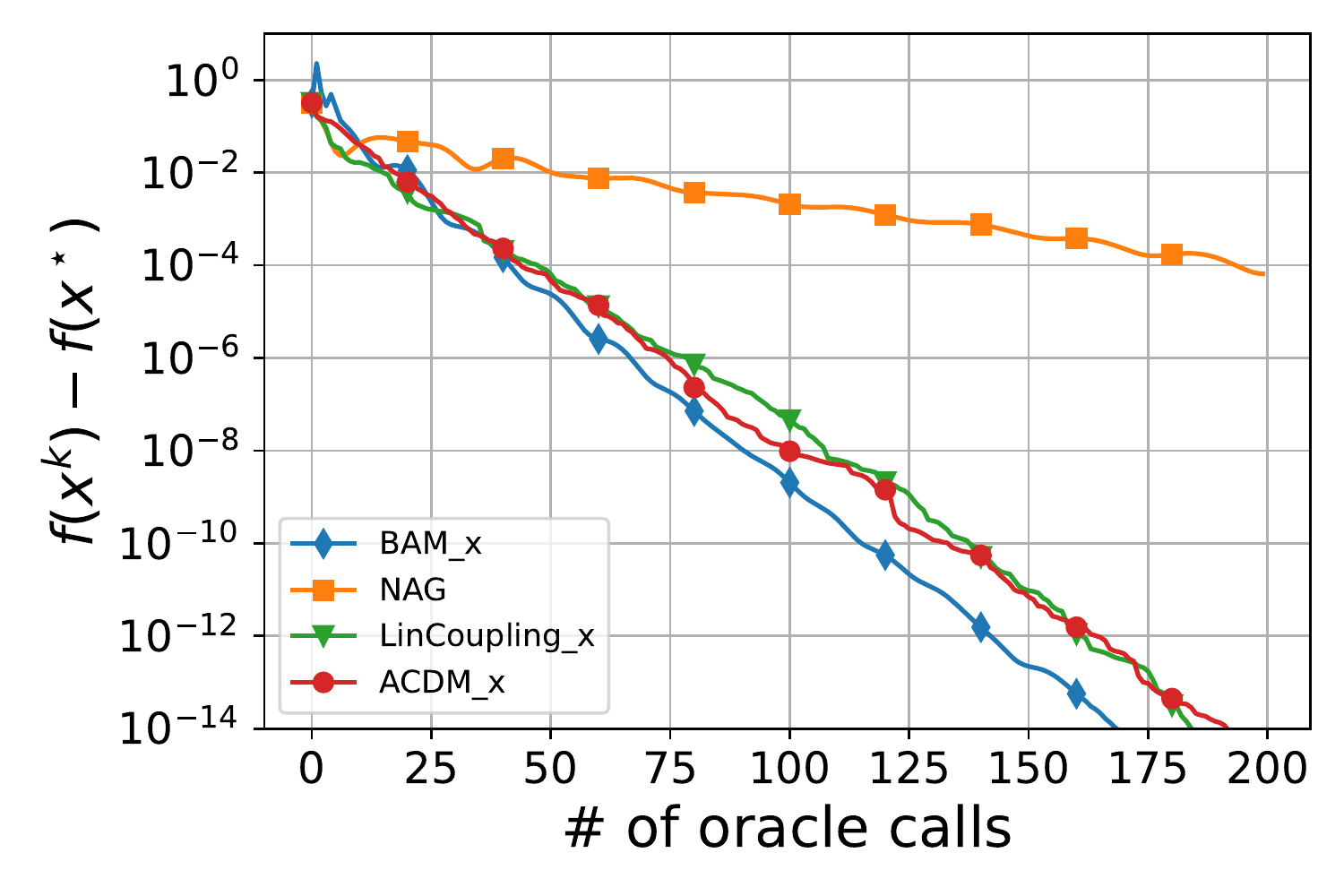}&	
		\includegraphics[width=0.31\linewidth]{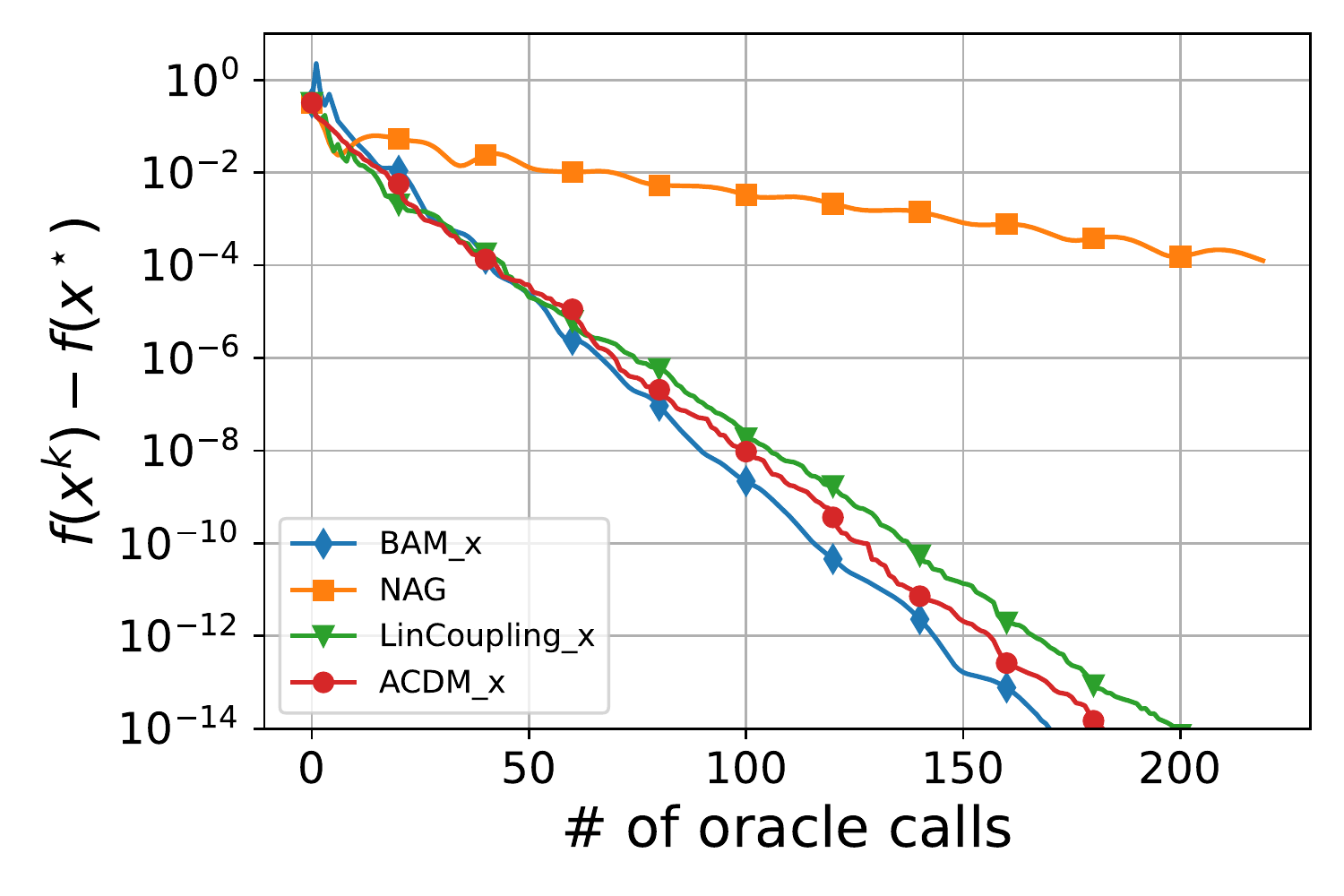}\\
  	    \includegraphics[width=0.31\linewidth]{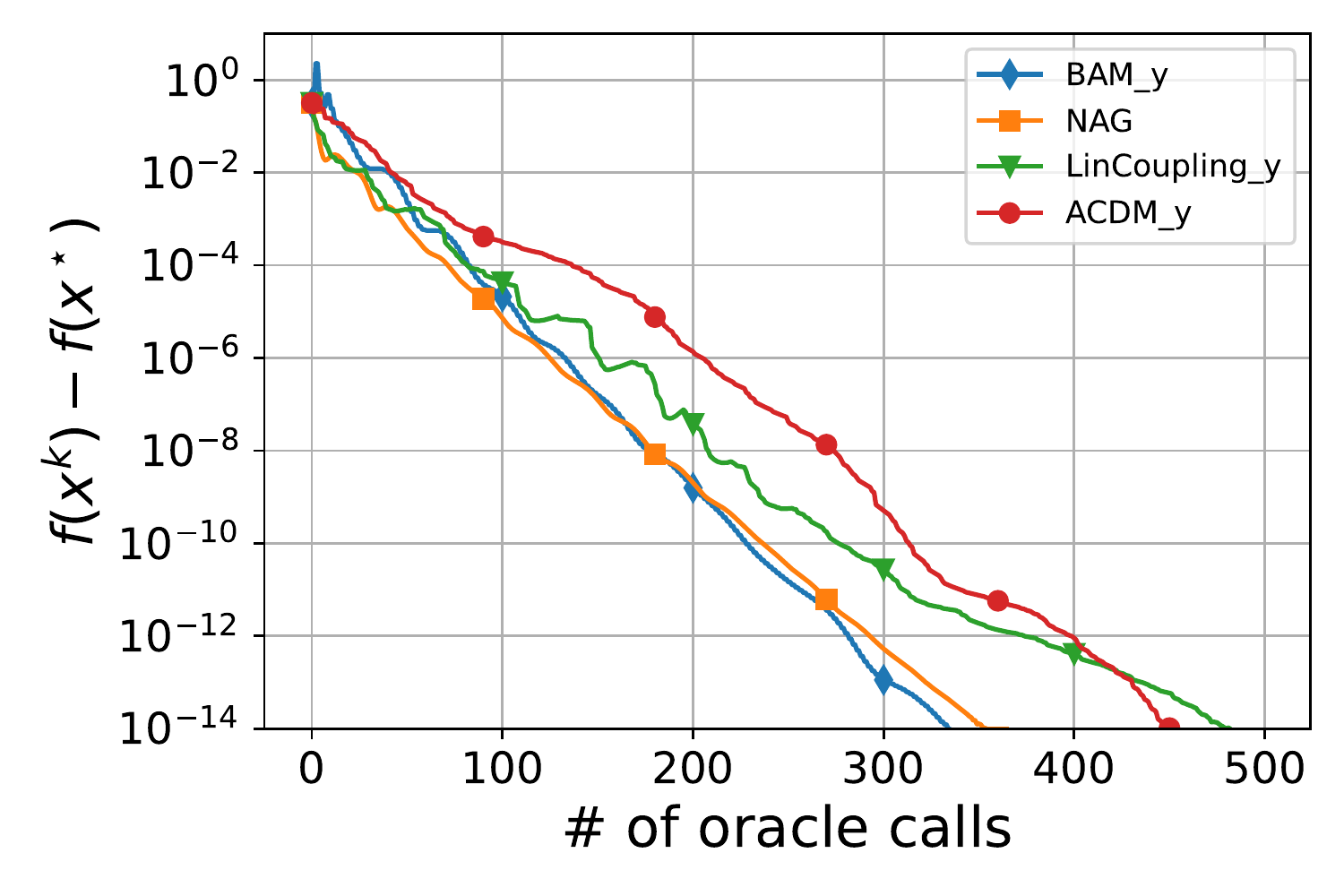}&	
        \includegraphics[width=0.31\linewidth]{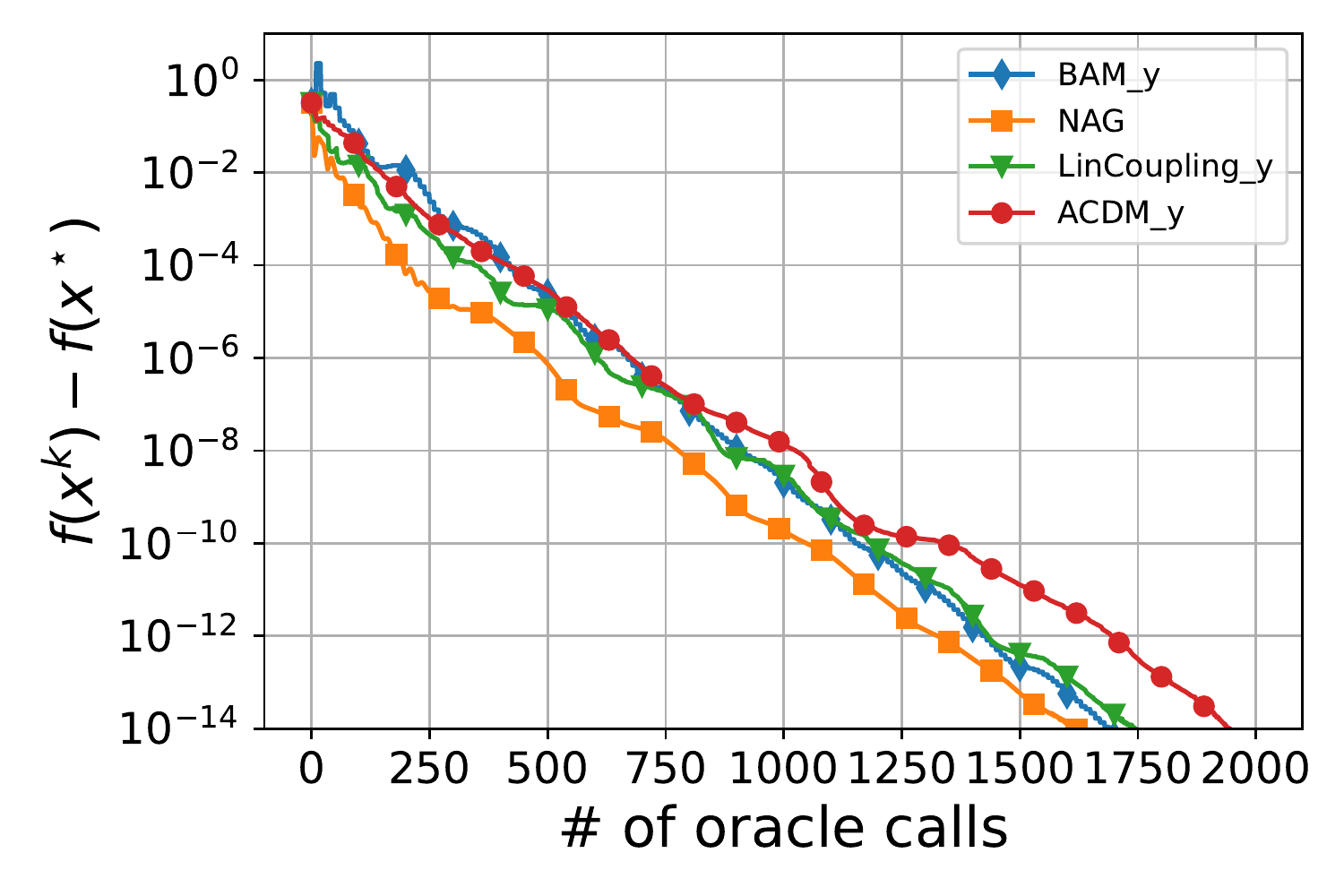}&	
		\includegraphics[width=0.31\linewidth]{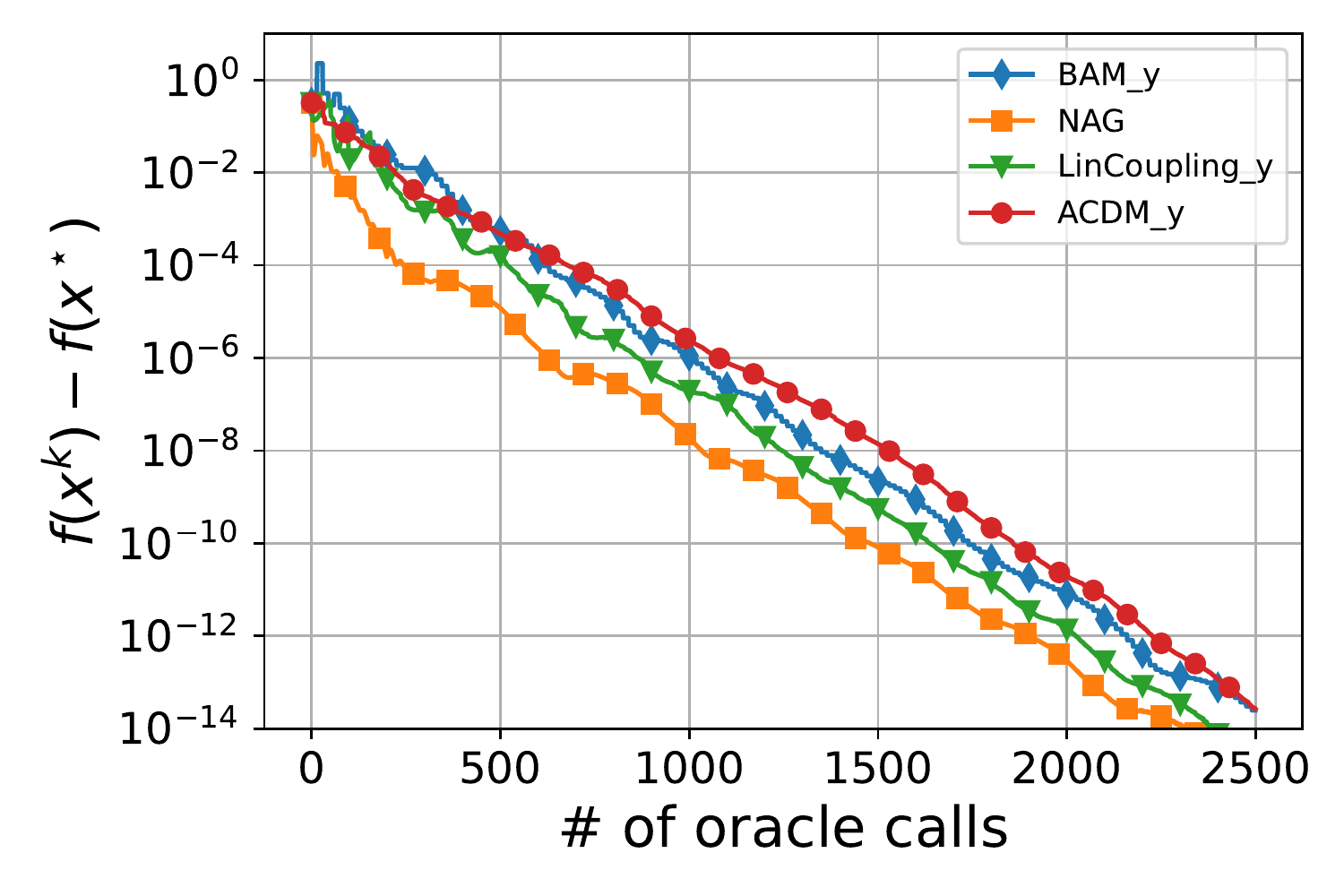}
	\end{tabular}
	\centering
	\caption{Comparison of Block Accelerated Method (BAM), Nesterov Accelerated Method (NAG), Accelerated Coordinate Descent Method (ACDM), and Linear Coupling method (LinCoupling) on logistic regression loss functions with two different $l_2$ regularizers. The first line represents the rate in terms of the $\nabla_x f(x,y)$ oracle calls, and the second one represents the rate in terms of the $\nabla_y f(x,y)$ oracle calls. We set $\mu_y = 0.002$ (left column), $\mu_y = 0.0001$ (middle column) and $\mu = 0.00005$ (right column). }
	\label{ris:image2}
\end{figure*}
\section{Experiments}
In all our experiments, we compare the new Block Accelerated Method (BAM) with Nesterov Accelerated Method (NAG)~\cite{nesterov1983method}, Accelerated Coordinate Descent Method (ACDM)~\cite{nesterov2017efficiency} and Linear Coupling method (LinCoupling)~\cite{allen2016even, gasnikov2015accelerated}. 
\subsection{Quadratic objectives}
In our experiments, we first consider quadratic functions:
$$
    f(z) = z^\top A z + b^\top z,
$$
where $z = (x,y)^\top$ is a joint vector of two blocks. The matrix spectrum is uniformly generated from $\mu_x$ to $L_x$ for the block $x$ of the matrix $A$ and from $\mu_y$ to $L_y$ for the block $y$ of matrix $A$. We set $\mu_x = \mu_y = 0.1$, $L_y = 50$ and we set dimensions of our blocks to be $d_x = 100$ and $d_y = 10$.  We vary the parameter $L_y$ to obtain different condition numbers $\kappa_y$, and we consider $\nabla_x f(x,y)$ and $\nabla_y f(x,y)$ oracle calls to compare several methods. 
\subsection{Logistic regression}
In our experiments, we also consider the logistic regression loss function with two $l_2$ regularizers for the click-prediction model:
\begin{align*}
f(x,y)&: = \frac{1}{n}\sum_{k=1}^n \log\left(1 + \exp\left(-\eta^k\langle\xi^k,(x,y)\rangle\right)\right)\\
    &+ \lambda_x \|x\|^2 + \lambda_y \|y\|^2.
\end{align*}
We used dataset "a1a" from LIBSVM collection~\cite{CC01a}. The smoothness constant of this dataset is estimated as $L = 1.567$. We set $d_x = 100$, $d_y = 19$ and $\mu_x = 0.01$. We vary the parameter $\mu_y$ to obtain different condition numbers $\kappa_y$. We also consider two oracle calls of $\nabla_x f(x,y)$ and $\nabla_y f(x,y)$.

\subsection{Results}
In our experiments, as seen on plots, the new method shows better performance in terms of $\nabla_x f(x,y)$ oracle calls for both objective functions and all condition numbers. Moreover, all accelerated coordinate methods outperform Nesterov Gradient Method significantly, which confirms theoretical bounds. In terms of $\nabla_y f(x,y)$ oracle calls, the new method shows approximately the same results as other accelerated coordinate methods and Nesterov Gradient Method. In case of expensive $\nabla_x f(x,y)$ oracle calls, the new method can be useful. Moreover, the new method can be generalized to distributed and federated settings, which means that this method has practical perspectives. 

\section{Discussion}
In this paper, we consider a convex optimization problem with a min-min structure 
$$\min_{x,y} f(x,y).$$
Assuming that $f$ is $L$-smooth and $\mu_x$-strongly convex in $x$, $\mu_y$-strongly convex in $y$ we propose new Algorithm BAM that required $\cO\left(\sqrt{L/\mu_x}\log\frac{1}{\epsilon} \right)$ calculations of $\nabla_x f$ and $\cO\left(\sqrt{L/\mu_y}\log\frac{1}{\epsilon} \right)$ calculations of $\nabla_y f$. The proposed in the paper approach allows different generalizations. For example, it can be generalized to mixed oracles \cite{gladin2021solving}: e.g., instead of $\nabla_y f$, only the value of $f$ is available. Another generalization is increasing the number of blocks (in this paper, we consider only two blocks, $x$, and $y$) for clarity. BAM can also be joined with many other tricks. For example, it can be joined with composite sliding \cite{lan2016gradient,kovalev2022optimal}, described at the very beginning of the introduction.


\nocite{langley00}

\bibliography{example_paper}
\bibliographystyle{icml2023}

\newpage
\appendix
\onecolumn
\section{Supplementary materials}

\begin{proof}[Proof of \Cref{lem:descent}]
	Using \cref{ass:L}, we get
	\begin{align*}
		f(\ol x^{k+1}, \ol y^{k+1})
		&\leq
		f(\ul x^{k}, \ol y^{k+1}) + \<\nabla_x f(\ul x^{k}, \ol y^{k+1}), \ol x^{k+1} - \ul x^k> + \frac{L_x}{2}\sqn{\ol x^{k+1} - \ul x^k}
		\\&=
		f(\ul x^{k}, \ol y^{k+1}) + \eta_x\alpha \left( \frac{\eta_x\alpha L_x}{2} - 1\right)\sqn{\nabla_x f(\ul x^{k}, \ol y^{k+1})}
		\\&\leq
		f(\ul x^{k}, \ol y^{k+1}) + \eta_x\alpha \left( \frac{1}{2} - 1\right)\sqn{\nabla_x f(\ul x^{k}, \ol y^{k+1})}
		\\&\leq
		f(\ul x^{k}, \ol y^{k+1}) - \frac{\eta_x\alpha}{2} \sqn{\nabla_x f(\ul x^{k}, \ol y^{k+1})}
	\end{align*}
\end{proof}

\begin{proof}[Proof of \Cref{thm:sliding}]
	\begin{align*}
		\eta_x^{-1}\cR_x^{k+1} + \eta_y^{-1}\cR_y^{k+1}
		&=
		\eta_x^{-1}\cR_x^{k} + \eta_y^{-1}\cR_y^{k}
		-\eta_x^{-1}\sqn{x^{k+1} - x^k}
		-\eta_y^{-1}\sqn{y^{k+1} - y^k}
		\\&
		+2\eta_x^{-1}\<x^{k+1} - x^k,x^{k+1} - x^*>
		+2\eta_y^{-1}\<y^{k+1} - y^k,y^{k+1} - y^*>
		\\&=
		\eta_x^{-1}\cR_x^{k} + \eta_y^{-1}\cR_y^{k}
		-\eta_x^{-1}\sqn{x^{k+1} - x^k}
		-\eta_y^{-1}\sqn{y^{k+1} - y^k}
		\\&
		+2\eta_x^{-1}\alpha\<\ul x_k - x^{k+1},x^{k+1} - x^*>
		+2\eta_y^{-1}\<\ol{y}^{k+1} - y^{k+1},y^{k+1} - y^*>
		\\&
		-2\<\nabla_x f(\ul x^k, \ol y^{k+1}),x^{k+1} - x^*>
		-2\<\nabla_y f(\ul x^k, \ol y^{k+1}),y^{k+1} - y^*>
		\\&=
		\eta_x^{-1}\cR_x^{k} + \eta_y^{-1}\cR_y^{k}
		-\eta_x^{-1}\sqn{x^{k+1} - x^k}
		-\eta_y^{-1}\sqn{y^{k+1} - y^k}
		\\&
		+\eta_x^{-1}\alpha\left(\sqn{\ul x^k - x^*} - \sqn{x^{k+1} - x^*} - \sqn{x^{k+1} - \ul x^k}\right)
		\\&
		+\eta_y^{-1}\alpha\left(\sqn{\ol y^{k+1} - y^*} - \sqn{y^{k+1} - y^*} - \sqn{y^{k+1} - \ol y^{k+1}}\right)
		\\&
		-2\<\nabla_x f(\ul x^k, \ol y^{k+1}),x^{k+1} - x^*>
		-2\<\nabla_y f(\ul x^k, \ol y^{k+1}),y^{k+1} - y^*>
		\\&\leq
		\eta_x^{-1}\cR_x^{k} + \eta_y^{-1}\cR_y^{k}
		-\eta_x^{-1}\sqn{x^{k+1} - x^k}
		-\eta_y^{-1}\sqn{y^{k+1} - y^k}
		\\&
		+\eta_x^{-1}\alpha\left(\sqn{\ul x^k - x^*} - \sqn{x^{k+1} - x^*} \right)
		\\&
		+\eta_y^{-1}\alpha\left(\sqn{\ol y^{k+1} - y^*} - \sqn{y^{k+1} - y^*} \right)
		\\&
		-2\<\nabla_x f(\ul x^k, \ol y^{k+1}),x^{k+1} - x^*>
		-2\<\nabla_y f(\ul x^k, \ol y^{k+1}),y^{k+1} - y^*>.
	\end{align*}
	This implies
	\begin{align*}
		(1+\alpha)\left(\eta_x^{-1}\cR_x^{k+1} + \eta_y^{-1}\cR_y^{k+1}\right)
		&\leq
		\eta_x^{-1}\cR_x^{k} + \eta_y^{-1}\cR_y^{k}
		-\eta_x^{-1}\sqn{x^{k+1} - x^k}
		-\eta_y^{-1}\sqn{y^{k+1} - y^k}
		\\&
		+\eta_x^{-1}\alpha\sqn{\ul x^k - x^*} 
		+\eta_y^{-1}\alpha\sqn{\ol y^{k+1} - y^*}
		\\&
		-2\<\nabla_x f(\ul x^k, \ol y^{k+1}),x^{k+1} - x^*>
		-2\<\nabla_y f(\ul x^k, \ol y^{k+1}),y^{k+1} - y^*>
		\\&=
		\eta_x^{-1}\cR_x^{k} + \eta_y^{-1}\cR_y^{k}
		-\eta_x^{-1}\sqn{x^{k+1} - x^k}
		-\eta_y^{-1}\sqn{y^{k+1} - y^k}
		\\&
		+\eta_x^{-1}\alpha\sqn{\ul x^k - x^*} 
		+\eta_y^{-1}\alpha\sqn{\ol y^{k+1} - y^*}
		\\&
		-2\<\nabla_x f(\ul x^k, \ol y^{k+1}),x^{k+1} - x^k>
		-2\<\nabla_x f(\ul x^k, \ol y^{k+1}),x^k - x^*>
		\\&
		-2\<\nabla_y f(\ul x^k, \ol y^{k+1}),y^{k+1} - y^k>
		-2\<\nabla_y f(\ul x^k, \ol y^{k+1}),y^k - y^*>.
	\end{align*}
	Using Young's inequality, we get
	\begin{align*}
		(1+\alpha)\left(\eta_x^{-1}\cR_x^{k+1} + \eta_y^{-1}\cR_y^{k+1}\right)
		&\leq
		\eta_x^{-1}\cR_x^{k} + \eta_y^{-1}\cR_y^{k}
		-\eta_x^{-1}\sqn{x^{k+1} - x^k}
		-\eta_y^{-1}\sqn{y^{k+1} - y^k}
		\\&
		+\eta_x^{-1}\alpha\sqn{\ul x^k - x^*} 
		+\eta_y^{-1}\alpha\sqn{\ol y^{k+1} - y^*}
		\\&
		+\eta_x^{-1}\sqn{x^{k+1} - x^k} + \eta_x\sqn{\nabla_x f(\ul x^k, \ol y^{k+1})}
		-2\<\nabla_x f(\ul x^k, \ol y^{k+1}),x^k - x^*>
		\\&
		+\eta_y^{-1}\sqn{y^{k+1} - y^k} + \eta_y\sqn{\nabla_y f(\ul x^k, \ol y^{k+1})}
		-2\<\nabla_y f(\ul x^k, \ol y^{k+1}),y^k - y^*>
		\\&=
		\eta_x^{-1}\cR_x^{k} + \eta_y^{-1}\cR_y^{k}
		+\eta_x^{-1}\alpha\sqn{\ul x^k - x^*} 
		+\eta_y^{-1}\alpha\sqn{\ol y^{k+1} - y^*}
		\\&
		+ \eta_x\sqn{\nabla_x f(\ul x^k, \ol y^{k+1})}
		+ \eta_y\sqn{\nabla_y f(\ul x^k, \ol y^{k+1})}
		\\&
		-2\<\nabla_x f(\ul x^k, \ol y^{k+1}),x^k - x^*>
		-2\<\nabla_y f(\ul x^k, \ol y^{k+1}),y^k - y^*>
		\\&=
		\eta_x^{-1}\cR_x^{k} + \eta_y^{-1}\cR_y^{k}
		+\eta_x^{-1}\alpha\sqn{\ul x^k - x^*} 
		+\eta_y^{-1}\alpha\sqn{\ol y^{k+1} - y^*}
		\\&
		+ \eta_x\sqn{\nabla_x f(\ul x^k, \ol y^{k+1})}
		+ \eta_y\sqn{\nabla_y f(\ul x^k, \ol y^{k+1})}
		\\&
		-2\<\nabla_x f(\ul x^k, \ol y^{k+1}),\ul x^k - x^*>
		-2\<\nabla_y f(\ul x^k, \ol y^{k+1}),\ol y^{k+1} - y^*>
		\\&
		-2\<\nabla_x f(\ul x^k, \ol y^{k+1}),x^k - \ul x^k>
		-2\<\nabla_y f(\ul x^k, \ol y^{k+1}),y^k - \ol y^{k+1}>
	\end{align*}
	Using \Cref{ass:mu}, we get
	\begin{align*}
		(1+\alpha)\left(\eta_x^{-1}\cR_x^{k+1} + \eta_y^{-1}\cR_y^{k+1}\right)
		&\leq
		\eta_x^{-1}\cR_x^{k} + \eta_y^{-1}\cR_y^{k}
		+\eta_x^{-1}\alpha\sqn{\ul x^k - x^*} 
		+\eta_y^{-1}\alpha\sqn{\ol y^{k+1} - y^*}
		\\&
		+ \eta_x\sqn{\nabla_x f(\ul x^k, \ol y^{k+1})}
		+ \eta_y\sqn{\nabla_y f(\ul x^k, \ol y^{k+1})}
		\\&
		+2\left(f(x^*, y^*) - f(\ul x^k, \ol y^{k+1})\right) - \mu_x \sqn{\ul x^k - x^*} - \mu_y \sqn{\ol y^{k+1} - y^*}
		\\&
		-2\<\nabla_x f(\ul x^k, \ol y^{k+1}),x^k - \ul x^k>
		-2\<\nabla_y f(\ul x^k, \ol y^{k+1}),y^k - \ol y^{k+1}>
		\\&=
		\eta_x^{-1}\cR_x^{k} + \eta_y^{-1}\cR_y^{k}
		+(\eta_x^{-1}\alpha - \mu_x)\sqn{\ul x^k - x^*} 
		+(\eta_y^{-1}\alpha - \mu_y)\sqn{\ol y^{k+1} - y^*}
		\\&
		+ \eta_x\sqn{\nabla_x f(\ul x^k, \ol y^{k+1})}
		+ \eta_y\sqn{\nabla_y f(\ul x^k, \ol y^{k+1})}
		+2\left(f(x^*, y^*) - f(\ul x^k, \ol y^{k+1})\right) 
		\\&
		-2\<\nabla_x f(\ul x^k, \ol y^{k+1}),x^k - \ul x^k>
		-2\<\nabla_y f(\ul x^k, \ol y^{k+1}),y^k - \ol y^{k+1}>
		\\&=
		\eta_x^{-1}\cR_x^{k} + \eta_y^{-1}\cR_y^{k}
		+(\eta_x^{-1}\alpha - \mu_x)\sqn{\ul x^k - x^*} 
		+(\eta_y^{-1}\alpha - \mu_y)\sqn{\ol y^{k+1} - y^*}
		\\&
		+ \eta_x\sqn{\nabla_x f(\ul x^k, \ol y^{k+1})}
		+ \eta_y\sqn{\nabla_y f(\ul x^k, \ol y^{k+1})}
		+2\left(f(x^*, y^*) - f(\ul x^k, \ol y^{k+1})\right) 
		\\&
		+\frac{2(1-\alpha)}{\alpha}\<\nabla_x f(\ul x^k, \ol y^{k+1}),\ol x^k - \ul x^k>
		\\&
		-\<\nabla_y f(\ul x^k, \ol y^{k+1}),\frac{2}{\alpha}(\ul y^k - (1-\alpha)\ol y^k)- \ol y^{k+1}>
		\\&=
		\eta_x^{-1}\cR_x^{k} + \eta_y^{-1}\cR_y^{k}
		+(\eta_x^{-1}\alpha - \mu_x)\sqn{\ul x^k - x^*} 
		+(\eta_y^{-1}\alpha - \mu_y)\sqn{\ol y^{k+1} - y^*}
		\\&
		+ \eta_x\sqn{\nabla_x f(\ul x^k, \ol y^{k+1})}
		+ \eta_y\sqn{\nabla_y f(\ul x^k, \ol y^{k+1})}
		+2\left(f(x^*, y^*) - f(\ul x^k, \ol y^{k+1})\right) 
		\\&
		+\frac{2(1-\alpha)}{\alpha}\<\nabla_x f(\ul x^k, \ol y^{k+1}),\ol x^k - \ul x^k>
		\\&
		+\frac{2}{\alpha}\<\nabla_y f(\ul x^k, \ol y^{k+1}),\ol y^{k+1}-\ul y^k + (1-\alpha)\ol y^k -(1-\alpha)\ol y^{k+1}>
		\\&=
		\eta_x^{-1}\cR_x^{k} + \eta_y^{-1}\cR_y^{k}
		+(\eta_x^{-1}\alpha - \mu_x)\sqn{\ul x^k - x^*} 
		+(\eta_y^{-1}\alpha - \mu_y)\sqn{\ol y^{k+1} - y^*}
		\\&
		+ \eta_x\sqn{\nabla_x f(\ul x^k, \ol y^{k+1})}
		+ \eta_y\sqn{\nabla_y f(\ul x^k, \ol y^{k+1})}
		+2\left(f(x^*, y^*) - f(\ul x^k, \ol y^{k+1})\right) 
		\\&
		+\frac{2(1-\alpha)}{\alpha}\left(\<\nabla_x f(\ul x^k, \ol y^{k+1}),\ol x^k - \ul x^k> + \<\nabla_y f(\ul x^k, \ol y^{k+1}),\ol y^k -\ol y^{k+1}>\right)
		\\&
		+\frac{2}{\alpha}\<\nabla_y f(\ul x^k, \ol y^{k+1}),\ol y^{k+1}-\ul y^k>.
	\end{align*}
	Using convexity of $f(x,y)$, we get
	\begin{align*}
		(1+\alpha)\left(\eta_x^{-1}\cR_x^{k+1} + \eta_y^{-1}\cR_y^{k+1}\right)
		&\leq
		\eta_x^{-1}\cR_x^{k} + \eta_y^{-1}\cR_y^{k}
		+(\eta_x^{-1}\alpha - \mu_x)\sqn{\ul x^k - x^*} 
		+(\eta_y^{-1}\alpha - \mu_y)\sqn{\ol y^{k+1} - y^*}
		\\&
		+ \eta_x\sqn{\nabla_x f(\ul x^k, \ol y^{k+1})}
		+ \eta_y\sqn{\nabla_y f(\ul x^k, \ol y^{k+1})}
		+2\left(f(x^*, y^*) - f(\ul x^k, \ol y^{k+1})\right) 
		\\&
		+\frac{2(1-\alpha)}{\alpha}\left(f(\ol x^k, \ol y^k) - f(\ul x^k, \ol y^{k+1})\right)
		+\frac{2}{\alpha}\<\nabla_y f(\ul x^k, \ol y^{k+1}),\ol y^{k+1}-\ul y^k>.
		\\&=
		\eta_x^{-1}\cR_x^{k} + \eta_y^{-1}\cR_y^{k}
		+(\eta_x^{-1}\alpha - \mu_x)\sqn{\ul x^k - x^*} 
		+(\eta_y^{-1}\alpha - \mu_y)\sqn{\ol y^{k+1} - y^*}
		\\&
		+ \eta_x\sqn{\nabla_x f(\ul x^k, \ol y^{k+1})}
		+ \eta_y\sqn{\nabla_y f(\ul x^k, \ol y^{k+1})}
		\\&
		+2 f(x^*, y^*)
		+\frac{2(1-\alpha)}{\alpha}f(\ol x^k, \ol y^k)
		-\frac{2}{\alpha}f(\ul x^k, \ol y^{k+1})
		\\&
		+\frac{2}{\alpha}\<\nabla_y f(\ul x^k, \ol y^{k+1}),\ol y^{k+1}-\ul y^k>.
	\end{align*}
	Using \Cref{lem:descent}, we get
	\begin{align*}
		(1+\alpha)\left(\eta_x^{-1}\cR_x^{k+1} + \eta_y^{-1}\cR_y^{k+1}\right)
		&\leq
		\eta_x^{-1}\cR_x^{k} + \eta_y^{-1}\cR_y^{k}
		+(\eta_x^{-1}\alpha - \mu_x)\sqn{\ul x^k - x^*} 
		+(\eta_y^{-1}\alpha - \mu_y)\sqn{\ol y^{k+1} - y^*}
		\\&
		+ \eta_x\sqn{\nabla_x f(\ul x^k, \ol y^{k+1})}
		+ \eta_y\sqn{\nabla_y f(\ul x^k, \ol y^{k+1})}
		\\&
		+2 f(x^*, y^*)
		+\frac{2(1-\alpha)}{\alpha}f(\ol x^k, \ol y^k)
		\\&
		-\frac{2}{\alpha}\left( f(\ol x^{k+1}, \ol y^{k+1}) + \frac{\eta_x\alpha}{2} \sqn{\nabla_x f(\ul x^{k}, \ol y^{k+1})}\right)
		\\&
		+\frac{2}{\alpha}\<\nabla_y f(\ul x^k, \ol y^{k+1}),\ol y^{k+1}-\ul y^k>
		\\&=
		\eta_x^{-1}\cR_x^{k} + \eta_y^{-1}\cR_y^{k}
		+(\eta_x^{-1}\alpha - \mu_x)\sqn{\ul x^k - x^*} 
		+(\eta_y^{-1}\alpha - \mu_y)\sqn{\ol y^{k+1} - y^*}
		\\&
		+\frac{2(1-\alpha)}{\alpha}\left(f(\ol x^k, \ol y^k) - f(x^*,y^*)\right)
		-\frac{2}{\alpha}\left( f(\ol x^{k+1}, \ol y^{k+1}) - f(x^*,y^*)\right)
		\\&
		+ \eta_y\sqn{\nabla_y f(\ul x^k, \ol y^{k+1})}
		+\frac{2}{\alpha}\<\nabla_y f(\ul x^k, \ol y^{k+1}),\ol y^{k+1}-\ul y^k>
		\\&=
		\eta_x^{-1}\cR_x^{k} + \eta_y^{-1}\cR_y^{k}
		+(\eta_x^{-1}\alpha - \mu_x)\sqn{\ul x^k - x^*} 
		+(\eta_y^{-1}\alpha - \mu_y)\sqn{\ol y^{k+1} - y^*}
		\\&
		+\frac{2(1-\alpha)}{\alpha}\left(f(\ol x^k, \ol y^k) - f(x^*,y^*)\right)
		-\frac{2}{\alpha}\left( f(\ol x^{k+1}, \ol y^{k+1}) - f(x^*,y^*)\right)
		\\&
		+ \eta_y\sqn{\nabla_y f(\ul x^k, \ol y^{k+1})}
		+2\eta_y\<\nabla_y f(\ul x^k, \ol y^{k+1}),(\eta_y\alpha)^{-1}(\ol y^{k+1}-\ul y^k)>
		\\&=
		\eta_x^{-1}\cR_x^{k} + \eta_y^{-1}\cR_y^{k}
		+(\eta_x^{-1}\alpha - \mu_x)\sqn{\ul x^k - x^*} 
		+(\eta_y^{-1}\alpha - \mu_y)\sqn{\ol y^{k+1} - y^*}
		\\&
		+\frac{2(1-\alpha)}{\alpha}\left(f(\ol x^k, \ol y^k) - f(x^*,y^*)\right)
		-\frac{2}{\alpha}\left( f(\ol x^{k+1}, \ol y^{k+1}) - f(x^*,y^*)\right)
		\\&
		+ \eta_y\sqn{\nabla_y f(\ul x^k, \ol y^{k+1})}
		+\eta_y\sqn{\nabla_y f(\ul x^k, \ol y^{k+1}) + (\eta_y\alpha)^{-1}(\ol y^{k+1}-\ul y^k)}
		\\&
		-\eta_y\sqn{\nabla_y f(\ul x^k, \ol y^{k+1})}
		-\eta_y^{-1}\alpha^{-2}\sqn{\ol y^{k+1} - \ul y^k}
		\\&=
		\eta_x^{-1}\cR_x^{k} + \eta_y^{-1}\cR_y^{k}
		+(\eta_x^{-1}\alpha - \mu_x)\sqn{\ul x^k - x^*} 
		+(\eta_y^{-1}\alpha - \mu_y)\sqn{\ol y^{k+1} - y^*}
		\\&
		+\frac{2(1-\alpha)}{\alpha}\left(f(\ol x^k, \ol y^k) - f(x^*,y^*)\right)
		-\frac{2}{\alpha}\left( f(\ol x^{k+1}, \ol y^{k+1}) - f(x^*,y^*)\right)
		\\&
		+\eta_y\left(\sqn{\nabla_y f(\ul x^k, \ol y^{k+1}) + (\eta_y\alpha)^{-1}(\ol y^{k+1}-\ul y^k)} - (\eta_y\alpha)^{-2}\sqn{\ol y^{k+1} - \ul y^k}\right).
	\end{align*}
	Using inequality \eqref{eq:ms}, we get
	\begin{align*}
		(1+\alpha)\left(\eta_x^{-1}\cR_x^{k+1} + \eta_y^{-1}\cR_y^{k+1}\right)
		&\leq
		\eta_x^{-1}\cR_x^{k} + \eta_y^{-1}\cR_y^{k}
		+(\eta_x^{-1}\alpha - \mu_x)\sqn{\ul x^k - x^*} 
		+(\eta_y^{-1}\alpha - \mu_y)\sqn{\ol y^{k+1} - y^*}
		\\&
		+\frac{2(1-\alpha)}{\alpha}\left(f(\ol x^k, \ol y^k) - f(x^*,y^*)\right)
		-\frac{2}{\alpha}\left( f(\ol x^{k+1}, \ol y^{k+1}) - f(x^*,y^*)\right).
	\end{align*}
	Using the choice of parameters $\eta_x,\eta_y,\alpha$, we get
	\begin{align*}
		(1+\alpha)\left(\eta_x^{-1}\cR_x^{k+1} + \eta_y^{-1}\cR_y^{k+1}\right)
		&\leq
		\eta_x^{-1}\cR_x^{k} + \eta_y^{-1}\cR_y^{k}
		+\frac{2(1-\alpha)}{\alpha}\left(f(\ol x^k, \ol y^k) - f(x^*,y^*)\right)
		\\&
		-\frac{2}{\alpha}\left( f(\ol x^{k+1}, \ol y^{k+1}) - f(x^*,y^*)\right).
	\end{align*}
	After rearranging, we get
	\begin{align*}
		\Psi^{k+1}
		&\leq
		\eta_x^{-1}\cR_x^{k} + \eta_y^{-1}\cR_y^{k}
		+\frac{2(1-\alpha)}{\alpha}\left(f(\ol x^k, \ol y^k) - f(x^*,y^*)\right)
		\\&\leq
		\eta_x^{-1}\cR_x^{k} + \eta_y^{-1}\cR_y^{k}
		+\frac{2(1+\alpha)^{-1}}{\alpha}\left(f(\ol x^k, \ol y^k) - f(x^*,y^*)\right)
		\\&=
		(1+\alpha)^{-1} \Psi^k.
	\end{align*}

\end{proof}

\begin{proof}[Proof of \Cref{cor:inner}]
	Using inequality~\eqref{eq:inner} and \eqref{eq:T}, we get
	\begin{align*}
		\norm{\nabla A^k(\ol y^{k+1})}
		&\leq
		\frac{(\eta_y\alpha)^{-1}}{2}\norm{\ul y^k- \argmin_{y \in \sY} A^k(y)}
		\\&\leq
		\frac{(\eta_y\alpha)^{-1}}{2}\norm{\ol y^{k+1} - \ul y^k}
		+\frac{(\eta_y\alpha)^{-1}}{2}\norm{\ol y^{k+1} - \argmin_{y \in \sY} A^k(y)}.
	\end{align*}
	Function $A^k(y)$ is $(\eta_y\alpha)^{-1}$-strongly convex which implies
	\begin{equation}
		(\eta_y\alpha)^{-1}\norm{\ol y^{k+1} - \argmin_{y \in \sY} A^k(y)} \leq \norm{\nabla A^k(\ol y^{k+1}) - \nabla A^k(\argmin_{y \in \sY} A^k(y))} = \norm{\nabla A^k(\ol y^{k+1})}.
	\end{equation}
	Hence, 
	\begin{align*}
		\norm{\nabla A^k(\ol y^{k+1})}
		&\leq
		\frac{(\eta_y\alpha)^{-1}}{2}\norm{\ol y^{k+1} - \ul y^k}
		+\frac{1}{2}\norm{\nabla A^k(\ol y^{k+1})}.
	\end{align*}
	Rearranging concludes the proof.
\end{proof}


\end{document}